\documentclass[12pt]{article}

\title{Set-valued expectiles for ordered data analysis}

\author{
Ha Thi Khanh Linh\footnote{Free University of Bolzano-Bozen, \href{mailto:thikhanhlinh.ha@unibz.it}{thikhanhlinh.ha@unibz.it}} \;
Andreas H. Hamel\footnote{Free University of Bolzano-Bozen, \href{mailto:andreas.hamel@unibz.it}{andreas.hamel@unibz.it}}
}
\date{{\small \today}}

\usepackage{array} 
\usepackage{verbatim} 
\usepackage{graphicx,float} 

\usepackage{caption} 
\usepackage{amsmath, amssymb, enumerate, color, stackrel}

\newtheorem{theorem}{Theorem}
\newtheorem{corollary}[theorem]{Corollary}
\newtheorem{remark}[theorem]{Remark}
\newtheorem{lemma}{Lemma}[subsection]
\newtheorem{definition}[theorem]{Definition}
\newtheorem{proposition}[theorem]{Proposition}
\newtheorem{example}[theorem]{Example}

\numberwithin{equation}{section}  
\numberwithin{figure}{section}    
\numberwithin{table}{section}     
\numberwithin{theorem}{section}

\newcommand{\abs}[1]{\ensuremath{\left| #1 \right|}}
\newcommand{\cb}[1]{\ensuremath{ \left\{ #1 \right\} }}

\newcommand{\bs}{\backslash}

\newcommand{\pend}{ \hfill $\square$ \medskip}

\renewcommand{\P}{\ensuremath{\mathcal{P}}}

\newcommand{\F}{\ensuremath{\mathcal{F}}}
\newcommand{\G}{\ensuremath{\mathcal{G}}}

\newcommand{\R}{\mathrm{I\negthinspace R}}

\newcommand{\N}{\mathrm{I\negthinspace N}}

\newcommand{\E}{\mathbb{E}}

\DeclareMathOperator{\esssup}{ess\,sup}
\DeclareMathOperator{\essinf}{ess\,inf}

\newcommand{\gr}{{\rm graph \,}}

\newcommand{\cl}{{\rm cl \,}}

\newcommand{\co}{{\rm co \,}}

\newcommand{\Int}{{\rm int\,}}

\newcommand{\One}{\mathrm{1\negthickspace I}}

\usepackage[
colorlinks=true, linkcolor=blue, citecolor=blue, urlcolor=blue, filecolor=blue
auskommentieren und unten pdfborder aktivieren]{hyperref}

\definecolor{color0}{gray}{.50}
\definecolor{color1}{rgb}{0,.2,.8}
\definecolor{color2}{rgb}{1,.2,0}
\definecolor{color3}{rgb}{.2,.7,.6}

\begin{document}

\maketitle

\begin{abstract}
Recently defined expectile regions capture the idea of centrality with respect to a multivariate distribution, but fail to describe the tail behavior while it is not at all clear what should be understood by a tail of a multivariate distribution. Therefore, cone expectile sets are introduced which take into account a vector preorder for the multi-dimensional data points. This provides a way of describing and clustering a multivariate distribution/data cloud with respect to an order relation. Fundamental properties of cone expectiles including dual representations of both expectile regions and cone expectile sets are established. It is shown that set-valued sublinear risk measures can be constructed from cone expectile sets in the same way as in the univariate case. Inverse functions of cone expectiles are defined which should be considered as rank functions rather than depth functions. Finally, expectile orders for random vectors are introduced and characterized via expectile rank functions.
\end{abstract}



\section{Introduction}

Expectile regression and expectiles were introduced in \cite{NeweyPowell87Economet} and gained popularity later on, see, for example, \cite{DeRossiHarvey09JE, Eilers13SM, HolzmannKlar16EJS} and in particular the literature review in \cite{DaouiaPaindaveine19R} with many further references. They also provide a new class of financial risk measures of the type introduced in \cite{ArtznerEtAl99MF}, which has been pointed out in \cite{BelliniEtAl14IME} along with dual representation results. Therein, it was also shown that expectile risk measures are basically the only sublinear cash-additive risk measures which are elicitable. In \cite{BelliniDiBernardino17EJF, ZaevskiNedeltchev23IRFA}, further applications of expectiles to risk management are discussed.

Expectiles were extended from the univariate to the multivariate case in \cite{DiGiorgiMcNeil16CSDA, CascosOchoa21JMVA} in form of set-valued functions although already Eilers  \cite{Eilers10Proc} gave a geometric construction for the bivariate case. In \cite{DiGiorgiMcNeil16CSDA, CascosOchoa21JMVA}, one can find the definition of expectile regions as another instance of depth regions whose development started with Tukey \cite{Tukey75PICM}. Alternative approaches can be found at different levels of generality in \cite{MaumeDeschampsRulliereSaid17DM, HerrmannHofertMailhot18SAJ, DaouiaPaindaveine19R}. 

When passing from the univariate to the multivariate case, the definition of order statistics and quantiles like the median becomes a conceptual obstacle because it is not a priori clear which order relation should replace the "natural" $\leq$-relation for real numbers. For this reason, central regions became surrogates for order statistics and quantiles in multivariate statistics. While this was quite successful with respect to characterizing multivariate distributions, it does not really help to answer questions like the following: What are the x\% best data points in a cloud of points in $\R^d$? What is the fraction of them which is better or worse than the median--however the latter is defined? Which data points are in the upper or lower tail of the distribution? The difficulty is that one cannot say which data points are better or worse than others. However, an order relation is often present in applications since decision makers have an understanding which (multi-dimensional) data points are better than others, and they often look for the best alternatives (or the least risky ones etc.). Therefore, new concepts are proposed which start with a vector order for the values of $\R^d$-valued random variables which is not complete in general, i.e., there may exist non-comparable points. The idea is similar to the one in \cite{HamelKostner18JMVA}, but different types of results are possible due to the convexity properties of expectiles which are not shared by quantiles. In particular, the sublinearity of set-valued expectiles allows a dual representation which also yields a dual description of the expectile region from \cite{DiGiorgiMcNeil16CSDA, CascosOchoa21JMVA} as a special case. In turn, this leads to new methods for computing cone expectiles, expectile regions and even univariate expectiles for finite data sets which are based on linear programing techniques.

In the next section, cone expectile regions are introduced and basic features established. A dual representation is given in Section \ref{Sec_DualRep}. This representation opens a straight path to the expectile regions from \cite{DiGiorgiMcNeil16CSDA, CascosOchoa21JMVA}: it is shown that the cone expectile sets, the expectile region discussed in \cite{CascosOchoa21JMVA} and families of scalar expectiles all are based on the same dual information. In Section \ref{SecExpectileRM}, set-valued expectile risk measures are introduced and the differences to the construction in \cite{CascosOchoa21JMVA} is explained. Section \ref{SecExpectileRank} includes the definition of cone expectile rank functions along with their basic properties. Moreover, new order relations for data points and stochastic dominance order random vectors related to expectiles are introduced which can be applied, e.g., in Multi-Criteria Decision Making and for the analysis of ordered data.

\section{Definition and basic properties}\label{Section_Definition}

In this paper, random variables $X \colon \Omega \to \R^d$ over a probability space $(\Omega, \F, \Pr)$ are considered. It is assumed that there is a vector preorder  for their values present, i.e., a reflexive and transitive relation on $\R^d$, which is compatible with addition and multiplication with non-negative numbers.

Such vector preorders are generated by non-trivial closed convex cones $C  \subseteq \R^d$, i.e., $\varnothing \neq C \neq \R^d$ via
\[
x \leq_C y \quad \Leftrightarrow \quad y - x \in C.
\]
In particular, $C = \{0\}$ furnishes the special case considered in \cite{CascosOchoa21JMVA}. In this case, points in $\R^d$ can only be compared if they are equal which is of course a very poor order relation. The case $C = \R^d_+$ provides the "natural" example of component-wise order. On the other end of the spectrum, one can use closed halfspaces $C = H^+(w) = \{z \in \R^d \mid w^\top  z \geq 0\}$ which generate complete order relations. A vector order is a vector preorder which is additionally antisymmetric which in turn is the case if $C \cap (-C) = \{0\}$, i.e., $C$ is pointed.

The presence of such a vector preorder is the feature that sets the developments in this paper apart from earlier contributions like \cite{Eilers10Proc, DiGiorgiMcNeil16CSDA, MaumeDeschampsRulliereSaid17DM, HerrmannHofertMailhot18SAJ, DaouiaPaindaveine19R, CascosOchoa21JMVA} and it also leads to new concepts and perspectives: it might be appropriate to consider the present paper (like \cite{HamelKostner18JMVA, HamelKostner22CSDA}) as a contribution to the analysis of ordered data. 

The set $C^+ = \{w \in \R^d \mid \forall z \in C \colon w^\top  z \geq 0\}$ is called the (positive) dual of the cone $C$ (sometimes also the polar cone). As usual, $L^1_d := L^1_d(\Omega, \F, \Pr)$ denotes the Banach space of integrable $\R^d$-valued random variables and  $\E[X]$ the component-wise expected value of $X$.

\begin{definition}
\label{Def_Cone_Exp}
Let $X \in L^1_d$ be an $\R^d$-valued random variable and $0 < \alpha < 1$. The set
\begin{equation}
\label{EqUCE}
E^\alpha_{-C}(X) = \bigcap_{w \in C^+} \cb{z \in \R^d \mid w^\top z \leq e_\alpha(w^\top X)}
\end{equation}
is called the downward cone $\alpha$-expectile and the set
\begin{equation}
\label{EqLCE}
E^{1-\alpha}_C(X) = \bigcap_{w \in C^+} \cb{z \in \R^d \mid w^\top z \geq e_{1-\alpha}(w^\top X)}
\end{equation}
the upward cone $\alpha$-expectile of $X$.
\end{definition}

\begin{figure}[t]
  \centering
  \includegraphics[width=0.9\textwidth]{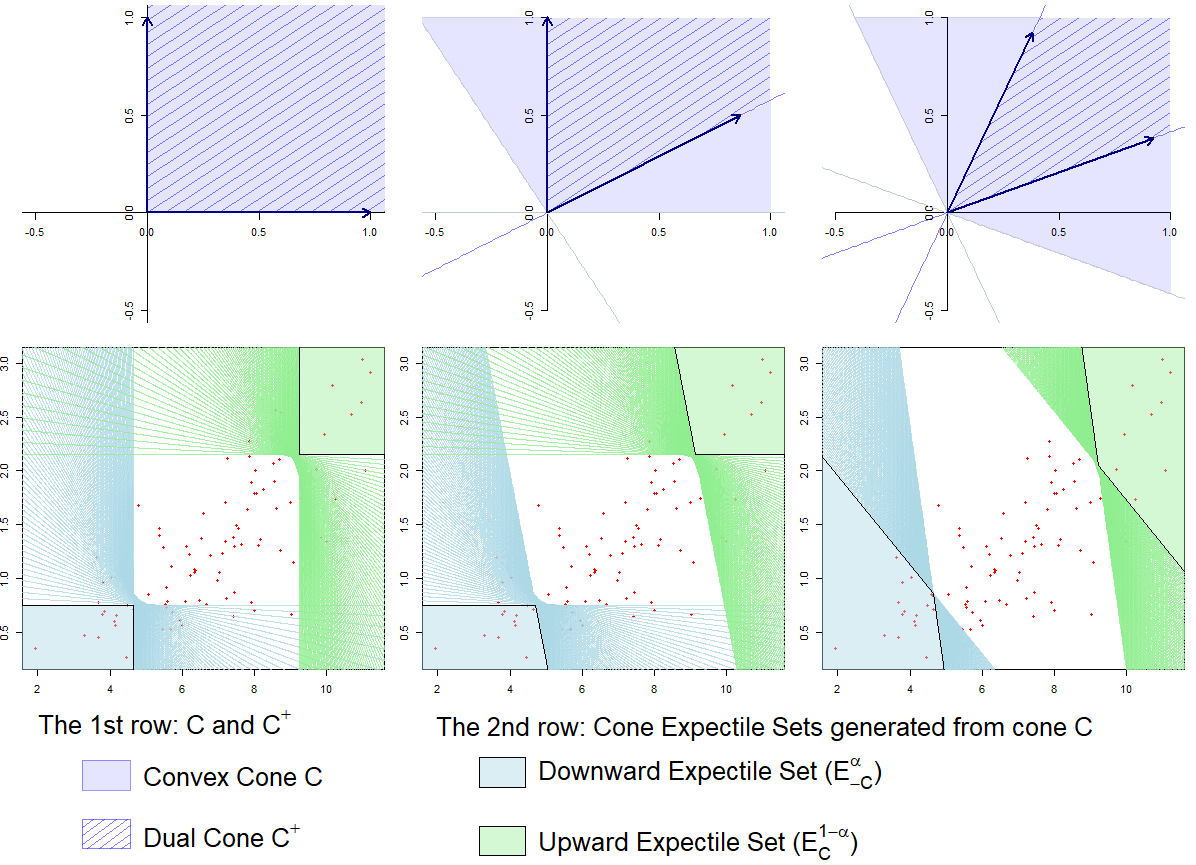}
  \caption{Examples of $E^\alpha_{-C}(X)$ and $E^{1-\alpha}_{C}(X)$ in $\R^2$ as intersections of halfspaces for different cones $C$ with $\alpha = 0.05$.}
  \label{EC_ConesDynamic}
\end{figure}

\begin{example}
\label{ExCones}
A sample of a bivariate variable simulated from the Gumbel copula in which the two marginal distributions are the normal distribution $N(7,4)$ and the gamma distribution $\Gamma(\alpha = 4, \beta =3)$ shows the shape of the cone expectiles for different cones $C$. This dataset, plotted as red points in Figure \ref{EC_ConesDynamic}, will also be used in the later examples. 

Figure \ref{EC_ConesDynamic} visualizes three different cones $C$ along with their duals in $\R^2$ (the first row) and how the choice of the cone $C$ affects the shape of downward and upward cone expectile sets (the second row). 

The computation of univariate sample expectiles for larger data sets in this paper was done in R using the package expectreg \cite{SobotkaEtAl14Github}.
\end{example}

\begin{remark}
A closed convex set $B^+ \subseteq C^+$ is said to be a base of the dual cone $C^+$ if for each $w \in C^+\bs\{0\}$ there are unique $s > 0$, $v \in B^+$ with $w = sv$. If $C^+$ has a base, then the intersections in \eqref{EqUCE} and \eqref{EqLCE} can be replaced by intersections over $B^+$ since the defining inequalities are positively homogeneous in $w$. The cone $C^+$ has a base if, and only if, $C$ has non-empty interior: in this case, the set $B^+ = \{w \in C^+ \mid w^\top  \bar z = 1\}$ is a base of $C^+$ where $\bar z \in \Int C$ can be chosen arbitrarily.
\end{remark}

In the special case $d=1$ and $C = \R_+$, the dual cone is $C^+ = C = \R_+$ with base $B^+ = \{1\}$ and one has
\begin{align*}
E^\alpha_{-\R_+}(X) & = \cb{z \in \R \mid z \leq e_\alpha(X)} = e_\alpha(X) -\R_+ \\
E^{1-\alpha}_{\R_+}(X) & = \cb{z \in \R \mid z \geq e_{1-\alpha}(X)} = e_{1-\alpha}(X) + \R_+.
\end{align*}

Another special case is $d=1$ and $C = \{0\}$ with dual cone $C^+ = \R$. In this case, it is enough to consider $w = 1$ and $w = -1$. One gets $E^\alpha_{-\{0\}}(X) = \cb{z \in \R \mid z \leq e_\alpha(X)} \cap \cb{z \in \R \mid -z \leq e_\alpha(-X) = - e_{1-\alpha}(-X)}$ as well as $E^{1-\alpha}_{\{0\}}(X) = \cb{z \in \R \mid z \geq e_{1-\alpha}(X)} \cap \cb{z \in \R \mid -z \geq e_{1-\alpha}(-X) = -e_\alpha(-X)}$. Thus, $E^\alpha_{-\{0\}}(X) = E^{1-\alpha}_{\{0\}}(X) = \varnothing$  for $0 < \alpha < \frac{1}{2}$ and $E^\alpha_{-\{0\}}(X) = E^{1-\alpha}_{\{0\}}(X) = [e_{1-\alpha}(X), e_\alpha(X)]$ for $\frac{1}{2} < \alpha < 1$. This interval is the singleton $\{\E[X]\}$ for $\alpha = \frac{1}{2}$. Of course, $E^\alpha_{-\{0\}}(X) = E^{1-\alpha}_{\{0\}}(X) = E^\alpha_{-\R_+}(X) \cap E^{1-\alpha}_{\R_+}(X)$.

\begin{remark}
\label{RemExpectationExpectile}
One has $e_\frac{1}{2}(X) = \E[X]$ if $d = 1$. This implies for $d > 1$
\[
E^\frac{1}{2}_{-C}(X) = \bigcap_{w \in C^+} \cb{z \in \R^d \mid w^\top  z \leq w^\top  \E[X]} = \E[X] - C.
\]
as well as
\[
E^{\frac{1}{2}}_C(X) = \bigcap_{w \in C^+} \cb{z \in \R^d \mid w^\top  z \geq w^\top  \E[X])} = \E[X] + C.
\]
\end{remark}

It is an important feature of $E^\alpha_{-C}$ and $E^{1-\alpha}_C$ that they map into well-defined complete lattices of subsets of $\R^d$, namely into
\begin{align*}
\G(\R^d, -C) & = \{D \subseteq \R^d \mid D = \cl\co(D - C)\} \; \text{and} \\
\G(\R^d, C) & = \{D \subseteq \R^d \mid D = \cl\co(D + C)\},
\end{align*}
respectively. It is well-known that $(\G(\R^d, -C), \subseteq)$ and $(\G(\R^d, C), \supseteq)$ are complete lattices (see \cite{HamelEtAl15Incoll}). The next results etablish a few elementary properties for the downward cone expectile.

\begin{proposition}
\label{PropUCEElementaryProps} The downward cone $\alpha$-expectiles satisfies:

(1) $E^\alpha_{-C}(X)$ is a closed convex set with $E^\alpha_{-C}(X) - C \subseteq E^\alpha_{-C}(X)$, i.e., $E^\alpha_{-C}(X) \in \G(\R^d, -C)$;

(2) If $A$ is an invertible $d \times d$-matrix and $b \in \R^d$, then $E^\alpha_{-AC}(AX + b) = AE^\alpha_{-C}(X) + b$ where $AD = \{Az \mid z \in D\}$ for a set $D \subseteq \R^d$;

(3) $X \leq_C Y$ $P$-a.s. implies $E^\alpha_{-C}(X) \subseteq E^\alpha_{-C}(Y)$;

(4) It holds $E^\alpha_{-C}(s X) = sE^\alpha_{-C}(X)$ for $s >0$, $X \in L^1_d$;

(5) It holds $E^\alpha_{-C}(X + Y) \supseteq E^\alpha_{-C}(X) + E^\alpha_{-C}(Y)$ for $0 < \alpha \leq \frac{1}{2}$, $X, Y \in L^1_d$;

(6) If $0 < \alpha \leq  \beta < 1$, then $E^\alpha_{-C}(X) \subseteq E^ \beta_{-C}(X)$;

(7) For $\alpha \in (0, 1)$, one has
\[
E^\alpha_{-C}(X) = \bigcap_{\beta \in (\alpha, 1)} E^\beta_{-C}(X).
\]
\end{proposition}

{\sc Proof.} (1) By definition, $E^\alpha_{-C}(X)$ is an intersection of closed halfspaces and therefore closed and convex. If $x \in E^\alpha_{-C}(X)$ and $z \in C$, then $w^\top  z \geq 0$ for all $w \in C^+$, hence $w^\top (x - z) \leq w^\top  x$ for all $w \in C^+$ which implies $x - z \in E^\alpha_{-C}(X)$. 

(2) The set $AC$ is a closed convex cone with dual $(AC)^+ = (A^\top )^{-1}C^+$ (see Lemma \ref{LemTransformDualCone} in the Appendix). Now,
\begin{align*}
z \in E^\alpha_{-AC}(AX) & \Leftrightarrow \forall v \in (AC)^+ \colon v^\top  z \leq e_\alpha(v^\top  (AX)) =  e_\alpha((A^\top  v)^\top  X) \\
	& \Leftrightarrow \forall w \in C^+ \colon w^\top  (A^{-1}z) = ((A^\top )^{-1}w)^\top  z \leq e_\alpha(w ^\top  X)) \\
	& \Leftrightarrow A^{-1}z \in E^\alpha_{-C}(X) \Leftrightarrow z \in A E^\alpha_{-C}(X)
\end{align*}
where the transformation $v = (A^\top )^{-1} w$ has been used. On the other hand, the translativity property of univariate expectiles \cite[Thm. 1 (iii)]{NeweyPowell87Economet}, \cite[Sect. 2.1, (iii)]{CascosOchoa21JMVA} gives $e_\alpha(w^\top  (X + b)) =  e_\alpha(w^\top  X) + w^\top  b$ for all $ w \in C^+$. This implies
\begin{align*}
& \{z \in \R^d \mid w^\top  z \leq e_\alpha(w^\top  (X + b)) = e_\alpha(w^\top  X) + w^\top  b\} \\
 & = \{z - b \in \R^d \mid w^\top  (z - b) \leq  e_\alpha(w^\top  X)\} + b = \{y \in \R^d \mid w^\top  y \leq  e_\alpha(w^\top  X)\} + b
\end{align*}
for all $ w \in C^+$. Taking the intersection over $w \in C^+$ gives $E^\alpha_{-C}(X + b) = E^\alpha_{-C}(X) + b$ which, together with the linear transformation property, proves (2).
 
(3) If $X \leq_C Y$ $P$-a.s., then $w^\top  X \leq w^\top  Y$  for all $w \in C^+$, hence $e_\alpha(w^\top  X) \leq e_\alpha(w^\top  Y)$ (see \cite[p. 3, property (v)]{CascosOchoa21JMVA}) which in turn implies $E^\alpha_{-C}(X) \subseteq E^\alpha(Y)$. 

(4) This follows directly from positive homogeneity of the univariate expectiles (see \cite[p. 3, property (iv)]{CascosOchoa21JMVA}). 

(5) The function $e_\alpha \colon L^1 \to \R$ is  superadditive for $0 < \alpha \leq \frac{1}{2}$ (see \cite[p. 3, property (viii)]{CascosOchoa21JMVA}), hence $e_\alpha(w^\top  X) + e_\alpha(w^\top  Y) \leq e_\alpha(w^\top  (X+Y))$ which in turn implies
\begin{multline*}
\cb{x \in \R^d \mid w^\top  x \leq e_\alpha(w^\top  (X+Y))}  \supseteq \cb{x \in \R^d \mid w^\top  x \leq e_\alpha(w^\top  X) + e_\alpha(w^\top  Y)} \\
	 = \cb{x \in \R^d \mid -w^\top  x \leq e_\alpha(w^\top  X)} + \cb{x \in \R^d \mid w^\top  x \leq e_\alpha(w^\top  Y)}
\end{multline*}
for all $w \in C^+$ where the last equality is a consequence of Lemma \ref{LemHalfSpaceSum}. Taking the intersection on both sides of the inclusion gives
\begin{equation}
\label{EqSupAdditive}
E^\alpha_{-C}(X + Y) \supseteq E^\alpha_{-C}(X) + E^\alpha_{-C}(Y)
\end{equation}
since the intersection of a sum always is a superset of the sum of the intersections. 

(6) This is a straightforward consequence of the parameter monotonicity of univariate expectiles: it holds $e_\alpha(w^\top  X) \leq e_\beta(w^\top  X)$ whenever $0 <  \alpha \leq \beta < 1$ for all $w \in C^+$ (see \cite[Proposition 5 (f)]{BelliniEtAl14IME}, \cite[p. 3, property (x)]{CascosOchoa21JMVA}).

(7) The inclusion $\subseteq$ is immediate from (6). Assume $z \in E^\beta_{-C}(X)$ for all $\beta \in (\alpha, 1)$. Then
\[
\forall w \in C^+\bs\{0\}, \forall \beta \in (\alpha, 1) \colon w^\top  z \leq e_\beta(w^\top  X).
\]
This implies $w^\top  z \leq e_\alpha(w^\top  X)$ for all $w \in C^+\bs\{0\}$ since $e_\beta(w^\top  X)$ is a continuous function of the parameter $\beta \in (0,1)$ (compare \cite[Proposition 1 (ii)]{HolzmannKlar16EJS}), hence $\supseteq$ is also true.
\pend

\begin{remark}
Properties (4) and (5) mean that the function $E^\alpha_{-C} \colon L^1_d \to \G(\R^d, -C)$ is positively homogeneous and superadditive for $0 < \alpha \leq \frac{1}{2}$ where the order relation $\subseteq$ in $\G(\R^d, -C)$ is understood as $\leq$.
\end{remark}

The corresponding result for upward cone $\alpha$-expectiles reads as follows. Its proof runs parallel to the one for Proposition \ref{PropUCEElementaryProps}.

\begin{proposition}
\label{PropLCEElementaryProps} The upward cone $\alpha$-expectiles satisfies:

(1) $E^{1-\alpha}_C(X)$ is a closed convex set with $E^{1-\alpha}_C(X) + C \subseteq E^{1-\alpha}_C(X)$, i.e., $E^{1-\alpha}_C(X) \in \G(\R^d, C)$;

(2) If $A$ is an invertible $d \times d$-matrix and $b \in \R^d$, then $E^{1-\alpha}_{AC}(AX + b) = AE^{1-\alpha}_C(X) + b$; 

(3) $X \leq_C Y$ $P$-a.s. implies $E^{1-\alpha}_C(X) \supseteq E^{1-\alpha}_C(Y)$;

(4) It holds $E^{1-\alpha}_C(s X) = sE^{1-\alpha}_C(X)$ for $s >0$, $X \in L^1_d$; 

(5) It holds $E^{1-\alpha}_C(X + Y) \supseteq E^{1-\alpha}_C(X) + E^{1-\alpha}_C(Y)$ for $X, Y \in L^1_d$, $0 < \alpha \leq \frac{1}{2}$;

(6) If $0 < \alpha \leq \beta < 1$, then $E^{1-\beta}_C(X) \supseteq E^{1-\alpha}_C(X)$. 

(7) For $\alpha \in (0, 1)$, one has
\[
E^{1-\alpha}_C(X) = \bigcap_{\beta \in (\alpha, 1)} E^{1-\beta}_C(X).
\]
\end{proposition}

Note that (5) now is subadditivity for the function $E^{1-\alpha}_C \colon  L^1_d \to \G(\R^d, C)$ since the order relation $\supseteq$ corresponds to $\leq$ in $\G(\R^d, C)$, not $\subseteq$. 

\begin{figure}
  \centering
  \includegraphics[width=0.6\textwidth]{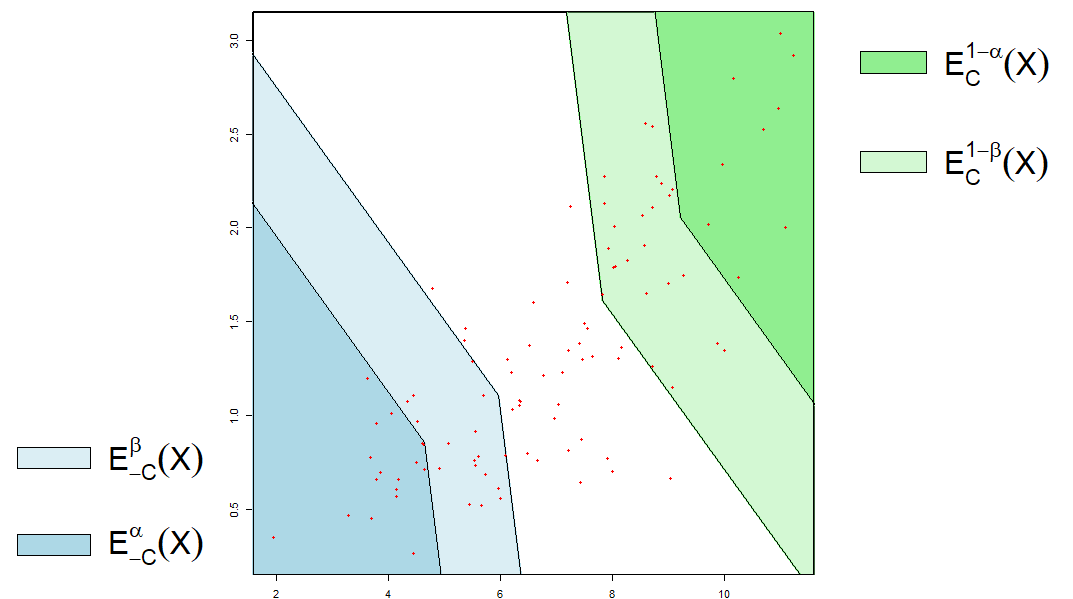}
  \caption{Comparing sets in lattices $\G(\R^2, C)$ and $\G(\R^2, -C)$.}
  \label{EC_Inclusion}
\end{figure}

Figure \ref{EC_Inclusion} provides an example of how the comparison of sets in the complete lattices $\G(\R^2, C)$ and $\G(\R^2, -C)$ works. For the downward cone expectile mapping into $\G(\R^2, -C)$, the set $E^{\alpha}_{-C}(X)$ is located "relatively lower" than the set $E^{\beta}_{-C}(X)$ for $\alpha < \beta$, i.e., $E^{\alpha}_{-C}(X) \subseteq E^{\beta}_{-C}(X)$. Therefore, $\subseteq$ stands for $\leq$ relation in lattice $\G(\R^2, -C)$. On the other hand, for the upward cone expectile mapping into $\G(\R^2, C)$ the set $E^{1-\beta}_{C}(X)$ is located "relatively lower" than the set $E^{1-\alpha}_{C}(X)$, i.e., $E^{1-\beta}_{C}(X) \supseteq E^{1-\alpha}_{C}(X)$. In this case, $\supseteq$ stands for $\leq$ relation in lattice $\G(\R^2, C)$. The cone $C$ here is strictly bigger than $\R^2_+$.

One should also realize that (4) and (5) in Proposition \ref{PropUCEElementaryProps} are equivalent to the fact that
\[
\gr E^\alpha = \cb{(X, x) \in L^1_d \times \R^d \mid x \in E^\alpha(X)}
\]
is a convex cone. Similarly, (4) and (5) in Proposition \ref{PropLCEElementaryProps} are equivalent to $\gr E^{1-\alpha}$ is a convex cone. Compare  \cite[Def. 4.1, Prop. 4.2]{HamelEtAl15Incoll} for a more general result of this type which covers the present situation. This cannot be concluded for the expectile region (see below for a definition) from positive homogeneity and (vii) in \cite[Sect. 3.1]{CascosOchoa21JMVA} for which the opposite inclusion in (5) is true. Moreover, while sublinearity in the previous two propositions is quite a natural feature, it turned out to be a serious issue in the geometric expectile approach: compare the involved constructions in \cite[Sect 3.2]{MaumeDeschampsRulliereSaid17DM} and \cite[Def. 4.4 ff]{HerrmannHofertMailhot18SAJ}.

The next result is another instance of a feature that is shared by upper/lower expectile sets on the one hand and univariate expectiles on the other hand. 

\begin{proposition}
\label{PropLUTransfer}
One has 
\begin{equation}
\label{EqLUTransfer}
\forall X \in L^1_d \colon E^{1-\alpha}_C(X) = -E^\alpha_{-C}(-X).
\end{equation}
\end{proposition}

{\sc Proof.} One has
\begin{align*}
-E^\alpha_{-C}(-X) & = \bigcap_{w \in C^+} \cb{-z \in \R^d \mid w^\top  z \leq e_\alpha(-w^\top  X) = -e_{1-\alpha}(w^\top  X)} \\
	& = \bigcap_{w \in C^+} \cb{y \in \R^d \mid -w^\top  y \leq  -e_{1-\alpha}(w^\top  X)} \\
	& = \bigcap_{w \in C^+} \cb{y \in \R^d \mid w^\top  y \geq e_{1-\alpha}(w^\top  X)} = E^{1-\alpha}_C(X)
\end{align*}
where in the first line \cite[p. 3, property (ii)]{CascosOchoa21JMVA} was used.
\pend

\begin{proposition}
\label{PropClosedGraph}
The graphs of the two functions $E^\alpha_{-C} \colon L^1_d \to \G(\R^d, -C)$ and $E^{1-\alpha}_C \colon L^1_d \to \G(\R^d, C)$ are closed in the product topology on $L^1_d \times \R^d$.
\end{proposition}

{\sc Proof.} Let $\{(X^n, z^n)\}_{n = 1, 2, \ldots} \subseteq \gr E^\alpha_{-C}$ be a sequence which converges to $(X, z) \in L^1_d \times \R^d$. Then $z^n \in E^\alpha_{-C}(X^n)$ for all $n = 1, 2, \ldots$ which is
\begin{equation}
\label{EqGraphConvergence}
\forall w \in C^+\bs\{0\}, \forall n \in \N  \colon w^\top  z^n \leq e_\alpha(w^\top  X^n).
\end{equation}
Theorem 10 in \cite{BelliniEtAl14IME} gives the Lipschitz continuity of the scalar $\alpha$-expectile for $0 < \alpha < 1$ with respect to the Wasserstein distance, i.e., one has for $\xi, \eta \in L^1_d$
\[
\abs{e_\alpha(\xi) - e_\alpha(\eta)} \leq \beta\inf\{\E\abs{\xi' - \eta'} \mid \xi' \sim P_\xi, \eta' \sim P_\eta\} \leq \beta \E\abs{\xi - \eta}
\]
where $\beta = \max\{\frac{\alpha}{1-\alpha}, \frac{1 - \alpha}{\alpha}\}$ and $\xi, \eta \in L^1$. Thus, $e_\alpha$ is continuous on $L^1$. This yields the convergence $e_\alpha(w^\top  X^n) \to e_\alpha(w^\top  X)$ for all $w \in C^+\bs\{0\}$. Hence, \eqref{EqGraphConvergence} implies
\[
\forall w \in C^+\bs\{0\}  \colon w^\top  z \leq e_\alpha(w^\top  X)
\]
which is $(X, z) \in \gr E^\alpha_{-C}$. The proof for $E^{1-\alpha}_C$ uses parallel arguments. 
\pend

\begin{figure}
  \centering
  \includegraphics[width=0.95\textwidth]{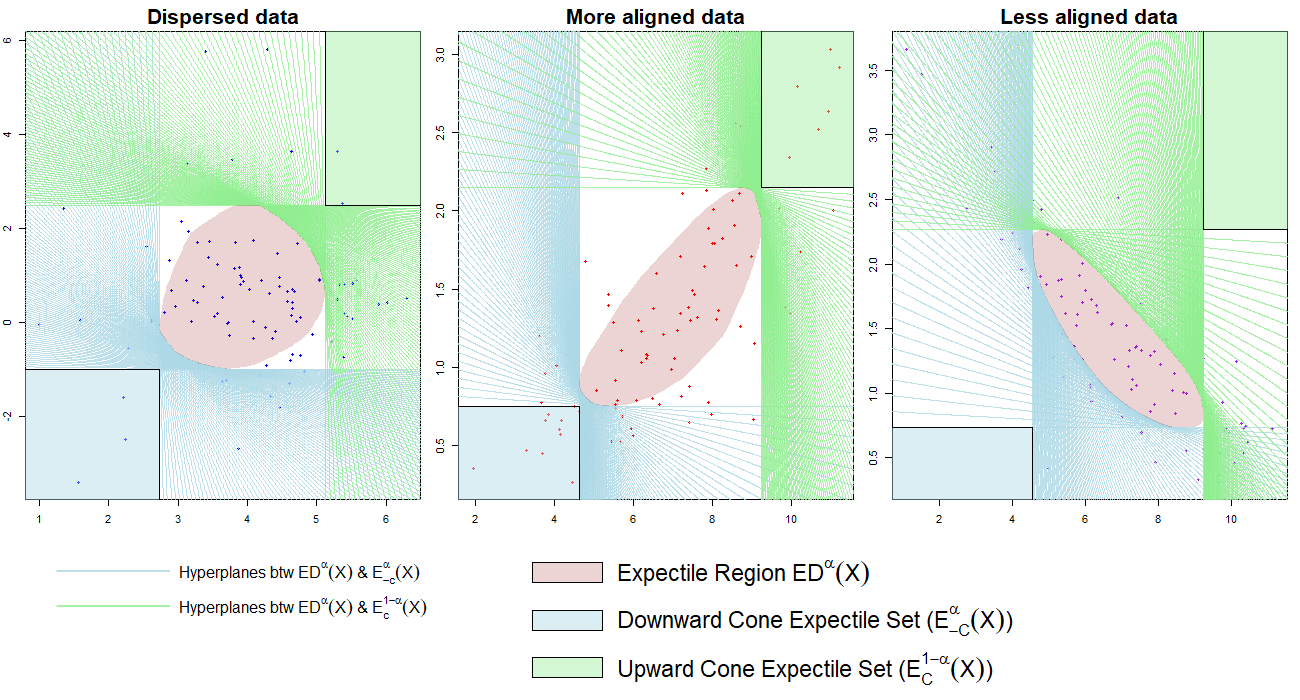}
  \caption{Expectile region and cone expectiles in $\R^2$ for $C=\R^2_+$ and $\alpha = 0.05$.}
   \label{ED_and_EC}
\end{figure}

This section is concluded with a clarification of the relationship between upward/downward cone expectiles and the expectile region from \cite{DiGiorgiMcNeil16CSDA, CascosOchoa21JMVA} which is defined as

\begin{equation} \label{EDdef}
ED^\alpha(X) = \bigcap_{w \in \mathbb S^{d-1}} \cb{z \in \R^d \mid w^\top  z \leq e_{1-\alpha}(w^\top  X)}
\end{equation}
for $0 < \alpha \leq \frac{1}{2}$ and $\mathbb S^{d-1} = \{w \in \R^d \mid w^\top  w = 1\}$. One has $z \in ED^\alpha(X)$ if, and only if,
\[
\forall w \in \R^d \colon w^\top  z \leq e_{1-\alpha}(w^\top  X).
\]
Thus, one also has
\[
\forall w \in \R^d \colon -w^\top  z \leq e_{1-\alpha}(-w^\top  X) = -e_\alpha(w^\top  X)
\]
which is
\[
\forall w \in \R^d \colon w^\top  z \geq e_\alpha(w^\top  X).
\]
In particular, it is always true for any $z \in ED^\alpha(X)$ that
\[
\forall w \in C^+ \colon e_\alpha(w^\top  X) \leq w^\top  z \leq e_{1-\alpha}(w^\top  X).
\]
This means that, for $0 < \alpha \leq \frac{1}{2}$, every  $w \in C^+\bs\{0\}$ produces a separating hyperplane for $E^\alpha_{-C}(X)$ and $ED^\alpha(X)$ on the one hand as well as $E^{1-\alpha}_{C}(X)$ and $ED^\alpha(X)$ on the other hand. The subsequent example, complemented by Figure \ref{ED_and_EC}, exemplifies this idea for various types of bivariate data sets.

\begin{example}
In Figure \ref{ED_and_EC}, three data samples are shown: the blue points in the left graph are simulated from two independent (univariate) random variables with normal distribution $N(4,1)$ and non-central Student's t-distribution $t(\nu=3, \delta =0.5)$.  The red points in the central graph are the same bivariate data as in Figures \ref{EC_ConesDynamic} and \ref{EC_Inclusion} above (see Example \ref{ExCones}). The purple points in the right graph are simulated from a rotated bivariate Gumbel Copula with the same parameters as the second one. It might be observed that the three samples can be considered as dispersed and more and less aligned, respectively, with the cone $C = \R^2_+$. In all three cases, the expectile region and the upward/downward cone expectiles are distinctly separated by hyperplanes generated from vectors $w \in C^+ = \R^2_+$.
\end{example}

From Figure \ref{ED_and_EC}, it can also be seen that the expectile region captures the central part of the multivariate distribution while downward and upward cone expectiles cover the lower and upper tails of the multivariate distribution, respectively, with respect to the order relation $\leq_C$ (the latter two may not even include a data point). One might say that the data points in $E^{\alpha}_{-C}(X)$ are the "very bad" points (the lower $\alpha$-expectile tail), the points in $E^{1-\alpha}_{C}(X)$ the "very good" ones (the upper $\alpha$-expectile tail); whereas one cannot categorize the remaining points for the given $\alpha$. It should be emphasized that "good" and "bad" are understood with maximization in view with respect to the order $\leq_{C}$ (vice versa for minimization). When increasing $\alpha$ from 0 to 1, the downward cone expectile set moves upward while the upward cone expectile set moves downward.

Both the expectile region $ED^\alpha(X)$ and the cone expectiles are set-valued functions of the random variable. Definition 1-2 in \cite{DaouiaPaindaveine19R} shares the same idea of set-valuedness but generalizes it to a family of M-quantiles which are defined as the minimizers of asymmetric loss functions. A different approach in terms of vector-valued expectiles for multivariate random variable can be found in \cite{HerrmannHofertMailhot18SAJ} and \cite{MaumeDeschampsRulliereSaid17DM}. Therein, the idea of elicitablility is used to produce (unique) minimizers of certain scoring functions on $\R^d$. While it might seem conceptually easier to consider a single vector as an expectile, this point hardly characterizes the multivariate distribution in a similar way as the scalar $\alpha$-expectile characterizes the univariate distribution because the distribution can be spread out in any direction in the $d$-dimensional space. 

In the quoted references, an ordering cone for the values of the random variable is not considered which also leads to  conceptual difficulties, for example, with respect to sublinearity as already discussed above (after Proposition \ref{PropLCEElementaryProps}). For this reason, cone expectiles together with the expectile region seem more appropriate since they can characterize the central as well as the tail behaviour with respect to a given vector order of a multivariate distribution.


\section{Dual Representation of Set-Valued Expectiles} \label{Sec_DualRep}

In this section, dual representation results are given for set-valued expectiles with $\alpha \in (0, \frac{1}{2})$. The univariate dual representation result from \cite[Proposition 8]{BelliniEtAl14IME} is used to represent set-valued expectiles as set-valued lower and upper expectations as discussed in \cite{HamelHeyde21Math}.

\begin{theorem}
\label{ThmDualRep}
Let $0 < \alpha \leq \frac{1}{2}$. Then
\[
e_\alpha(w^\top  X) = \inf\cb{\E^Q[w^\top  X] \mid Q \in \mathcal W(\alpha)} \; \text{as well as} \;
	e_{1-\alpha}(w^\top  X) = \sup\cb{\E^Q[w^\top  X] \mid Q \in \mathcal W(\alpha)} 
\]
where $\mathcal W(\alpha) = \cb{Q \in \mathcal M^P_1 \mid \frac{dQ}{dP} \in L^\infty, \, \esssup \frac{dQ}{dP} \leq \beta \essinf \frac{dQ}{dP}}$ with $\beta = \max\cb{\frac{\alpha}{1-\alpha}, \frac{1-\alpha}{\alpha}}$ and one has
\begin{align*}
E^\alpha_{-C}(X)  & =  \bigcap_{w \in C^+} \cb{z \in \R^d \mid w^\top  z \leq  \inf\cb{\E^Q[w^\top  X] \mid Q \in \mathcal W(\alpha)}} \\
	& = \bigcap_{w \in C^+} \bigcap_{Q \in \mathcal W(\alpha)} \cb{z \in \R^d \mid w^\top  z \leq \E^Q[w^\top  X]}
\end{align*}
as well as
\begin{align*}
E^{1-\alpha}_{C}(X)  & =  \bigcap_{w \in C^+} \cb{z \in \R^d \mid \sup\cb{\E^Q[w^\top  X] \mid Q \in \mathcal W(\alpha)} \leq w^\top  z} \\
	& = \bigcap_{w \in C^+} \bigcap_{Q \in \mathcal W(\alpha)} \cb{z \in \R^d \mid \E^Q[w^\top  X] \leq w^\top  z}.
\end{align*}
\end{theorem}

{\sc Proof.} The first two formulas follow from \cite[Proposition 8]{BelliniEtAl14IME} with a straighforward transformation of $L^\infty$- density functions into probability measures $Q$. The formulas for $E^\alpha_{-C}$ and $E^{1-\alpha}_{C}$ are now consequences of the definitions. \pend

\begin{remark}
\label{RemDualExpectile}
For fixed $w \in C^+\backslash\cb{0}$, $Q \in \mathcal{W}(\alpha)$, the function 
\[ 
X \mapsto E^+_{Q,w}(X) := \cb{z \in \R^d \mid w^\top \E^Q[X] \leq w^\top z}
\]
is an upper $(Q,w)$-expectation as defined in \cite[Sect. 4.7]{HamelHeyde21Math} since, of course, $\E^Q[w^\top X] = w^\top \E^Q[X]$ holds true where $\E^Q[X]$ is understood as the component-wise expectation of $X$ with respect to $Q$. The previous theorem establishes that the upward cone expectile is a set-valued upper expectation (and the corresponding downward cone expectile is a set-valued lower expectation). This again is a very striking parallelity to the univariate case in \cite[Proposition 8]{BelliniEtAl14IME}.
\end{remark}

The dual representations ask for the infimum and supremum of the linear function $w \mapsto w^\top \E^Q[X]$ over the set
\begin{equation}
\label{Set_EX}
W_\alpha(X) := \cb{\E^Q[X] \mid Q \in \mathcal W(\alpha)}.
\end{equation}

This raises the question what shape this set has and where it can be found. The following example provides experimental evidence.

\begin{example}
\label{ExEQSet} 
Let $\Omega = \{\omega_1, \omega_2, \omega_3 \}$, $p_i = P(\{\omega_i\})=\frac{1}{3}$, $i = 1,2,3$. This gives
\[
\begin{aligned}
\mathcal W(\alpha) &= \cb{(q_1, q_2, q_3) \in [0,1]^3 \mid q_1 + q_2 + q_3 = 1, \, \max\cb{\frac{q_1}{p_1}, \frac{q_2}{p_2}, \frac{q_3}{p_3}} \leq \beta \min\cb{\frac{q_1}{p_1}, \frac{q_2}{p_2}, \frac{q_3}{p_3}}} \\
&= \cb{(q_1, q_2, q_3) \in [0,1]^3 \mid q_1 + q_2 + q_3 = 1, \, \max\cb{q_1, q_2, q_3} \leq \beta \min\cb{q_1, q_2, q_3}}\\
\end{aligned}
\]
for $\beta = \max \cb{\frac{\alpha}{1-\alpha}, \frac{1-\alpha}{\alpha}}$. If $\alpha \in \left(0, \frac{1}{2}\right)$, then $\beta = \frac{1-\alpha}{\alpha}$ which is the case considered in this example.

Figure \ref{Set_Q} shows how the set of scenarios $\mathcal W(\alpha)$ looks for four different $\alpha$'s ($\alpha \in \cb{ 0.05, 0.25, 0.4, 0.499}$). In each case, the shape of the set $\mathcal W(\alpha)$ is filled by 10,000 (maroon) points $Q=\{q_1, q_2, q_3\}$ which will be used for computing the dual representations.

\begin{figure}
  \centering
  \includegraphics[width=\textwidth]{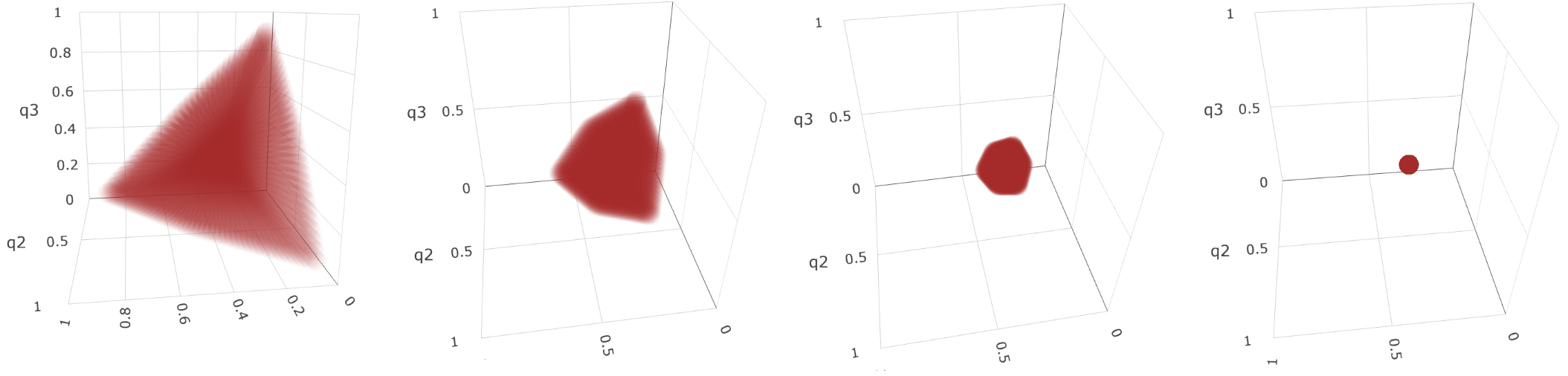}
  \vspace{-6pt}
  \caption{The set $\mathcal W(\alpha)$ of scenarios Q for 3 data points in Example \ref{ExEQSet} with $\alpha = 0.05, 0.25, 0.4, 0.499$ (from left to right).}
   \label{Set_Q}
\end{figure}

A discrete random variable $X \colon \Omega \to \R^2$  is given by: $\{x_1= X(\omega_1), x_2= X(\omega_2), x_3= X(\omega_3)\}$ $=\cb{\begin{bmatrix}5 \\ 2 \\ \end{bmatrix} , \begin{bmatrix}4 \\ -1 \\ \end{bmatrix},  \begin{bmatrix} 3 \\ 1 \\ \end{bmatrix} } \subset \R^2$.  
The computations in this simple example are done by brutal force for a fixed $\alpha=0.25$. The 10,000 sample points in the set $\mathcal W(\alpha = 0.25)$ are used to compute the dual representation of the cone expectile sets as per Theorem \ref{ThmDualRep}. Figure \ref{Dual_EC_graph} depicts the downward and upward cone expectile sets for 3 datapoints which are computed both from the definitions via \eqref{EqUCE}, \eqref{EqLCE} (the far left graph), and from the dual representation in Theorem \ref{ThmDualRep} (the far right graph). As expected, the definition and the dual representation give the same result.

\begin{figure}
  \centering
  \includegraphics[width=0.75\textwidth]{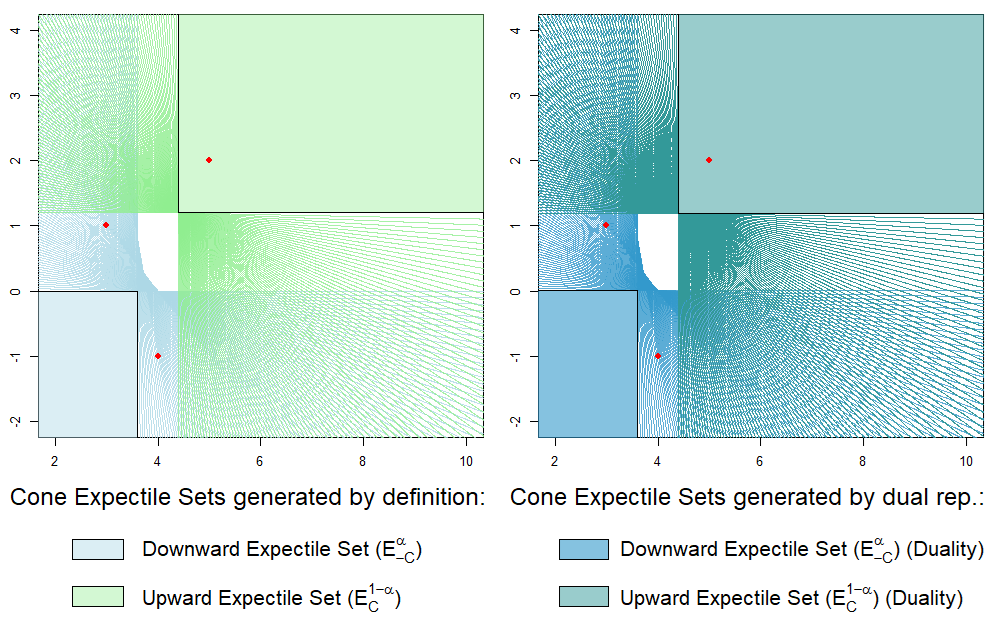}
  \vspace{-6pt}
  \caption{Downward and upward cone expectiles for 3 data points in Example \ref{ExEQSet}.}
   \label{Dual_EC_graph}
\end{figure}

The set $W_{\alpha}(X)$ is computed in a similar way: for each probability measure $Q$ (each maroon point in Figure \ref{Set_Q}, the second graph from left), $\E^Q[X]$ is computed. Then 10,000 points of $\E^Q[X]$ are plotted on the same plot of expectile region (the blue points in Figure \ref{ED_W(X)}). These blue points perfectly fill up the shape of the expectile region $ED^\alpha(X)$ defined in \cite{DiGiorgiMcNeil16CSDA, CascosOchoa21JMVA} (or in \eqref{EDdef} above). 

\begin{figure}
  \centering
  \includegraphics[width=0.75\textwidth]{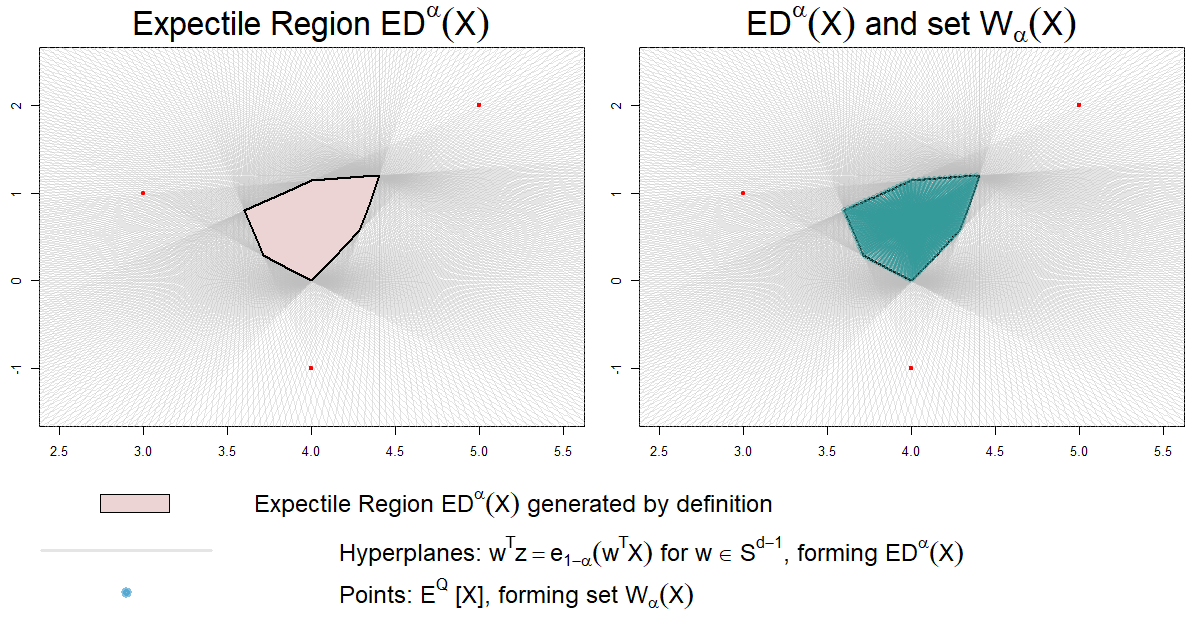}
  \vspace{-6pt}
  \caption{The set $W_{\alpha}(X)$ and the expectile region $ED^\alpha(X)$ for 3 data points  in Example \ref{ExEQSet} with $\alpha = 0.25$.}
   \label{ED_W(X)}
\end{figure}
\end{example}

As suggested by Example \ref{ExEQSet}, especially Figure \ref{ED_W(X)}, the following theorem establishes the coincidence of the set $W_\alpha(X)$ with the expectile depth region $ED^\alpha(X)$ defined in \cite{DiGiorgiMcNeil16CSDA, CascosOchoa21JMVA}. 

\begin{theorem}
\label{ThmDualExpectileDepthRegion}
For $X \in L^1_d$, one has
\[
\forall \alpha \in (0, \frac{1}{2} ]  \colon  W_\alpha(X) = ED^\alpha(X).
\]
\end{theorem}

{\sc Proof.} By definition,
\[
ED^\alpha(X) = \bigcap_{w \in \R^d\bs\{0\}} \cb{z \in \R^d \mid w^\top  z \leq e_{1-\alpha}(w^\top  X)}
\]
The dual representation of $e_{1-\alpha}(w^\top  X)$ provides
\[
z \in ED^\alpha(X) \quad \Leftrightarrow \quad \forall w \in \R^d \colon w^\top  z \leq \sup_{Q \in \mathcal W(\alpha)}\E^Q[w^\top  X].
\]
Now, the function
\[
w \mapsto \sup_{Q \in \mathcal W(\alpha)}w^\top \E^Q[X] = \sup_{y \in W_\alpha(X)} w^\top  y
\]
is the support function of the set $W_\alpha(X)$. According to a basic duality theorem (see \cite[5.83 Corollary]{AliprantisBorder06Book}) one has
\[
W_\alpha(X) = \bigcap_{w \in \R^d} \cb{z \in \R^d \mid w^\top  z \leq \sup_{y \in W_\alpha(X)} w^\top  y}.
\]
since $W_\alpha(X)$ is closed and convex (and compact). \pend

The previous theorem shows that the knowledge of the set $\mathcal W(\alpha)$ admits to reconstruct not only $E^\alpha_{-C}(X)$, $E^{1-\alpha}_{C}(X)$, but also $W_\alpha(X) = ED^\alpha(X)$. This is a dual way of constructing the expectile regions from \cite{CascosOchoa21JMVA} which also opens a path to use linear programming techniques for computing sample expectiles. Moreover, the set $\mathcal W(\alpha)$ is independent of $X$ and $C$ while $W_\alpha(X)$ depends on $X$ but still not on $C$. This can be seen as structuring the modelling process: if one chooses to use expectiles at some level $\alpha$, one gets set $\mathcal W(\alpha)$ of probability measures from the dual representation. Then, combined with the "pure" data $X$, one gets $W_\alpha(X)$--the analyst/decision maker is only interested in the shape of the distribution. Finally, if the analyst/decision maker has a preference for the data points, the cone $C$ enters the picture and one can additionally filter out the "(very) good" and the "(very) bad" data points with respect to the order $\leq_C$.

In \cite[Sec. 5]{CascosOchoa21JMVA}, an algorithm for computing the bivariate expectile region is proposed. However, it remains a challenge to compute expectile regions and lower/upper cone expectile sets for dimensions $d>2$. A general methodology is indicated in the following which draws on the dual representations in Theorem \ref{ThmDualRep} and \ref{ThmDualExpectileDepthRegion}. More specific and also more effective algorithms will be discussed elsewhere.

The main difficulty when computing expectile regions in $\R^d$ for $d \geq 2$ is to sample the vectors $w$ in all directions of the unit sphere $S^{d-1}$. This could be overcome by using the dual representation of the expectile region given in \eqref{Set_EX} and Theorem \ref{ThmDualExpectileDepthRegion} in case of finitely many data points, e.g., a sample of a distribution. In this case, $W_\alpha(X)$ will be a polytope in $\R^d$ which is--through the dual representation--the projection of another one in $\R^N$ where $N$ is the number of data points. This means it would be sufficient to compute the vertices of $W_\alpha(X) = ED^\alpha(X)$. The following result makes this idea more precise.

\begin{proposition}
\label{PropComputeED}
Let $\Omega = \cb{\omega_1, \omega_2, \ldots, \omega_N}$ be a finite set. Then $\mathcal W(\alpha)$ can be identified with a polytope which is a subset of the simplex $S = \cb{q \in \R^N_+ \mid \sum_{n=1}^N q_n = 1} \subseteq \R^N$.

If $X \colon \Omega \mapsto \R^d$ is a random variable, then $W_\alpha(X) = \cb{\E^Q[X] \mid Q \in \mathcal W(\alpha)} \subseteq  \R^d$ is a polytope. Moreover, if $z \in \R^d$ is a vertex of $W_\alpha(X)$, then there is a vertex $P$ of $\mathcal W(\alpha)$ such that $z = \E^P[X]$.
\end{proposition}

{\sc Proof.} The identification $Q \sim q \in \R^N$ via $q_n = Q(\omega_n)$ for $n \in \{, \ldots, N\}$ together with the fact that in the finite probability space setting the set $\mathcal W(\alpha)$ is defined by a finite set of linear equations and inequalities verifies the first claim. 

Next, assume $z = \E^P[X]$ is a vertex of $W_\alpha(X)$ and $P$ is not a vertex of $\mathcal{Q}$. Then $P$ can be written as a convex combination of some elements $Q_k \in \mathcal{Q}, k = 1, \ldots, K$ which are vertices of $\mathcal W(\alpha)$, i.e., $P= s_1 Q_1 + s_2 Q_2 +...+ s_K Q_K$ with $s_k \geq 0$ and $\sum s_k = 1$. One obtains $\E^P[X] = \E^{s_1 Q_1 +...+ s_K Q_K}[X] = s_1 \E^{ Q_1}[X] + ... + s_K \E^{ Q_K}[X]$. Since $\E^P[X]$ is a vertex of $W_\alpha(X)$, this implies $\E^P[X] = \E^{ Q_1}[X] = ... = \E^{ Q_K}[X]$, and each of the $Q_k$'s yields the vertex $z = \E^P[X]$. \pend

The brutal force idea would be to compute the vertices of $\mathcal W(\alpha)$ in $\R^N$, i.e., get a so called vertex representation, and then use it to compute the vertices of $W_\alpha(X)$. This is, however, not very efficient since $N$, i.e., the number of data points, is usually much larger than the dimension $d$. Therefore, it might be better to understand the representation $W_\alpha(X) = ED^\alpha(X)$ as a polyhedral projection problem and solve it either directly or as an equivalent multi-criteria linear programming problem. Efficient algorithms are available for such problems \cite{LoehneWeissing16MMOR, Weissing20MMOR}. 

For cone expectile sets with $0 < \alpha \leq \frac{1}{2}$, the direct computation from Definition \ref{Def_Cone_Exp} might be feasible since sampling $w \in C^+$ is no longer an issue if $C^+$ is polyhedral cone. The next theorem gives the key idea.

\begin{theorem}
\label{PolyhedralConeExp}
If $C^{+}$ is a polyhedral convex cone generated by $\cb{w^1, \ldots, w^M}$, i.e., $C^+ = \{w \in \R^d \mid w = \sum_{m = 1}^{M} s_m w^m, \, s_1, \ldots, s_M \geq 0\}$, then one has
\[
E^{\alpha}_{-C}(X)= \bigcap_{w^m \in \{w^1, \ldots, w^M\}} \cb{z \in \R^d \mid (w^m)^\top z \leq e_{\alpha}(w^m)^\top X)}
\]
for $0< \alpha \leq \frac{1}{2}$. A parallel formula holds for the downward cone expectile set $E^{1-\alpha}_{C}(X)$.
\end{theorem}

{\sc Proof.} See Appendix \ref{Prof_PolConeExp}. 
\pend

Theorem \ref{PolyhedralConeExp} states that, in the polyhedral case, cone expectile sets can be represented by intersections over a finite number of halfspaces instead of the infinitely many generated by all vectors $w \in C^+$. Therefore, the computation of cone expectile sets is simplified if the preference of the decision maker is given by a polyhedral cone $C$.


\section{Expectile risk measures for multivariate positions}
\label{SecExpectileRM}

In \cite{BelliniDiBernardino17EJF}, properties of the expectile risk measure defined by
\[
\varrho^{exp}_\alpha(X) = -e_\alpha(X)
\]
are studied; especially, this risk measure is monotone, cash additive and sublinear for $0 < \alpha \leq \frac{1}{2}$. Moreover, it is known that it is the only coherent risk measure among those which can be written as generalized quantiles \cite{BelliniEtAl14IME}. The corresponding result for set-valued expectile risk measure reads as follows. Compare \cite[Definition 6]{HamelHeyde21Math} and \cite[Definition 7.6]{HamelEtAl15Incoll} for the definition of set-valued risk measures.

\begin{theorem}
\label{ThmExpectRiskProperties}
If $C \supseteq \R^d_+$, then the function $R^{exp}_\alpha \colon L^1_d \to \P(\R^d, C)$ defined by $R^{exp}_\alpha(X) := -E^\alpha_{-C}(X) = E^{1-\alpha}_{C}(-X)$ with $0 < \alpha \leq \frac{1}{2}$ is a sublinear set-valued risk measure, i.e.,

(1) $R^{exp}_\alpha(0) \not\in \{\R^d, \varnothing\}$ (finiteness at zero),

(2) $R^{exp}_\alpha(X + z\One) = R^{exp}_\alpha(X) - z$ for all $X \in L^1_d$, for all $z \in \R^d$ (translativity),

(3) $X \leq_C Y$ $P$-a.s. implies $R^{exp}_\alpha(Y) \supseteq R^{exp}_\alpha(X)$ (monotonicity),

(4) $R^{exp}_\alpha(s X) = sR^{exp}_\alpha(X)$ for all $s > 0$ (positive homogeneity),

(5) $R^{exp}_\alpha(X + Y) \supseteq R^{exp}_\alpha(X) + R^{exp}_\alpha(Y)$ (subadditivity).

(6) $R^{exp}_\alpha$ has values in $\G(\R^d, C)$ and $\gr R^{exp}_\alpha \subseteq  L^1_d \times \R^d$ is closed in the product topology.
\end{theorem}

{\sc Proof.} One has
\begin{align*}
-E^\alpha_{-C}(X) = \bigcap_{w \in C^+} \cb{- x \in \R^d \mid w^\top  x \leq e_\alpha(w^\top  X)} 
 = \bigcap_{w \in C^+} \cb{z \in \R^d \mid w^\top  z \geq -e_\alpha(w^\top  X)}.
\end{align*}

(1) $R^{exp}_\alpha(0) = \bigcap_{w \in \R^d_+} \cb{z \in \R^d \mid w^\top  z \geq 0} = \R^d_+$ since $-e_\alpha(0) = 0$. (2) follows from Proposition \ref{PropUCEElementaryProps} with $A$ the $d \times d$ unit matrix. (3)-(5) also follow from Proposition \ref{PropUCEElementaryProps}, (3)-(5). (6) follows from Proposition \ref{PropUCEElementaryProps} (1) and Proposition \ref{PropClosedGraph}.
\pend

The elements of $R^{exp}_\alpha(X)$ are those deterministic portfolios which can be used as compensation for the risk of $X$ at initial time. The assumption $C \supseteq \R^d_+$ ensures in particular that one can add a non-negative portfolio to one that is already risk compensating and still get a risk compensating portfolio. Moreover, (3) applies in particular if $X - Y \in \R^d_+$ $P$-a.s.

For a univariate $X$, one has $\E[-X] = -e_{0.5}(X) \leq -e_\alpha(X) \leq -e_\beta(X)$ for $0 < \beta \leq \alpha \leq \frac{1}{2}$. This implies 
\[
R^{exp}_{0.5}(X) \supseteq R^{exp}_\alpha(X) \supseteq R^{exp}_\beta(X)
\]
for the multivariate case. This means that $\alpha = \frac{1}{2}$ provides the most optimistic risk measure, whereas smaller $\alpha$ means more risk aversity since there are potentially less options to compensate the risk. Note that this is right the opposite of the nesting property for expectile (and other depth) regions as stated in \cite[Sect. 3.1, (ii)]{CascosOchoa21JMVA}. Moreover, this also shows that the set $ED^\alpha(X) + \R^d_+$ is not appropriate as collection of risk compensating portfolios as implied in \cite[Sect. 3.1]{CascosOchoa21JMVA}: it always includes $\E[X]$ which already in the univariate case is not an appropriate risk compensating quantity. Moreover, the monotonicity property for the expectile region, namely $X \leq_{\R^d_+} Y$ $P$-a.s. implies $ED^\alpha(X) + \R^d_+ \supseteq ED^\alpha(Y) + \R^d_+$, would mean that there are less risk compensating portfolios available for the more risky $X$ compared to $Y$. The same point can be raised for the subadditivity property of expectile regions in \cite[Sect. 3.1 (vii)]{CascosOchoa21JMVA}: it would actually punish diversification. Overall, $X \mapsto ED^\alpha(X) + \R^d_+$ cannot be considered as an appropriate risk measure, but $R^{exp}_\alpha$ can which is also a direct generalization of the univariate expectile risk measure.

Finally note that the cone $C \subseteq \R^d$ can be understood as a solvency cone at initial time (see Section \ref{Section_Definition}), which also justifies the assumption $\R^d_+ \subseteq C$. It is an obvious question how such a cone (or even more general market models) can be taken into account. General procedures for this can be found in \cite{HamelHeydeRudloff11MAFE, HamelHeyde21Math}. For expectile risk measures, this will be discussed elsewhere.

The dual representation of the expectile risk measure is an immediate corollary of Theorem \ref{ThmDualRep}.

\begin{corollary}\label{CorExpectRiskDual}
For $X \in L^1_d$ and $0<\alpha \leq \frac{1}{2}$, one has 
\[
R^{exp}_{\alpha}(X) = \bigcap_{w \in C^+} \bigcap_{Q \in \mathcal{W}(\alpha)} \cb{z \in \R^d \mid \E^Q[w^\top (-X)] \leq w^\top  z}
\]
\end{corollary}

{\sc Proof.} This is a direct consequence of Theorem \ref{ThmDualRep} and the definition of $\R^{exp}_{\alpha}$. \pend

Remark \ref{RemDualExpectile} also yields that the expectile risk measure is a set-valued upper expectation--which is the core of dual representation results for coherent risk measures in general established in the seminal paper \cite{ArtznerEtAl99MF}. This, in turn, can be traced back to \cite[Sect. 10.2]{Huber81Book}.

\section{Cone expectile rank functions}
\label{SecExpectileRank}

Depth functions in multivariate statistics were invented with the idea of centrality: they quantify how close a given point is to the center of the distribution which in turn can be a single point or a (closed convex) set of points, the (most) central region. However, if one understands points as "good" if they dominate many other points with respect to the order $\leq_C$ generated by a cone $C$ and as "bad" if they are dominated by many other points, then "centrality" basically means "average behaviour" which corresponds to the fact that the expected value is an element of all expectile regions in the sense of \cite{CascosOchoa21JMVA}. 

Decision makers are usually not interested in average alternatives, but in "really good" ones with respect to their preferences where "good" depends on the goal (maximization or minimization). In such situations, centrality concepts and depth functions are of little use and should be complemented by functions which are monotone increasing with respect to the underlying order $\leq_C$. The following definition provides such functions generated by cone expectile sets. They relate to cone expectiles like the expectile depth function in \cite[Sect. 4]{CascosOchoa21JMVA} relates to expectile regions \cite[Sect. 3]{CascosOchoa21JMVA}.

\begin{definition}
\label{DefExpectileRank}
The functions $D_{-C}(\cdot; X), D_C(\cdot; X) \colon \R^d \to [0, 1]$ defined by
\begin{align*}
D_{-C}(z; X) & = \inf\cb{\alpha \in (0, 1) \mid z \in  E^\alpha_{-C}(X)} \\
D_C(z; X) & = \sup\cb{\beta \in (0, 1) \mid z \in  E^\beta_{C}(X)}
\end{align*}
are called downward and upward expectile rank function generated by $X$, respectively.
\end{definition}

For $d=1$ and $C=C^+=\R_+$, one has $D_{-C}(z; X) = D_C(z; X) = e^{-1}_{X}(z)$ which is a direct consequence of Definition \ref{DefExpectileRank}. Here, $e^{-1}_{X}(\cdot)$ is the inverse expectile function (see \cite[Sec. 2.1]{CascosOchoa21JMVA}). 

\begin{figure}[h]
  \centering
  \includegraphics[width=0.85\textwidth]{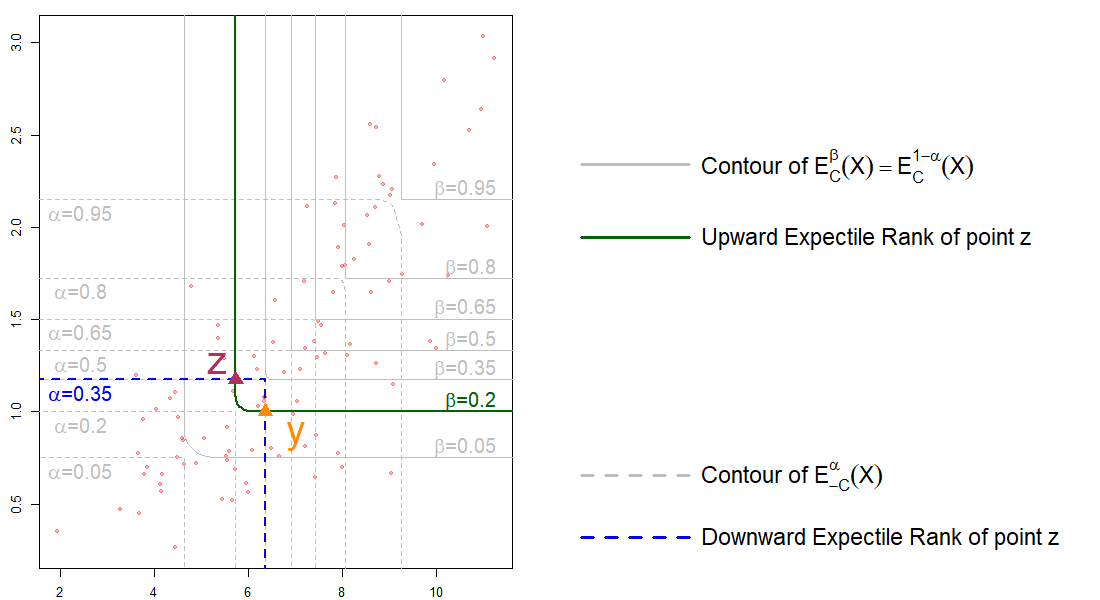}
  \caption{Downward and upward expectile ranks of points for a bivariate distribution.}
   \label{Cone_Exp_Rank}
\end{figure}

Figure \ref{Cone_Exp_Rank} shows an example of how to spot the downward and upward expectile ranks of a point in $\R^2$ generated by $X \in L^1_2$ for $C = \R^2_+$. The downward expectile rank of the "red" point $z$ is $D_{-\R^2_+}(z; X) = 0.35$ and its upward expectile rank is $D_{\R^2_+}(z; X) = 0.2$. 

The next result characterizes the cone expectile sets as the sublevel and superlevel sets of the downward and upward expectile rank function, respectively.

\begin{theorem}
\label{ThmExpectileRegionVsRank}
One has for $\alpha, \beta \in (0, 1)$
\begin{align*}
z \in E^\alpha_{-C}(X) \quad &\Leftrightarrow \quad D_{-C}(z; X) \leq \alpha, \\
z \in E^\beta_C(X) \quad & \Leftrightarrow \quad D_C(z; X) \geq \beta.
\end{align*}
\end{theorem}

{\sc Proof.} If $z \in E^\alpha_{-C}(X)$, then clearly $D_{-C}(z; X) \leq \alpha$. Vice versa, assume $D_{-C}(z; X) \leq \alpha$. The definition of the infimum and the monotonicity of the expectile set as in Proposition \ref{PropUCEElementaryProps} (6) imply $z \in E^{\alpha + \varepsilon}_{-C}(X)$ for all $\varepsilon > 0$ with $\alpha + \varepsilon < 1$. Hence, $z \in  E^\alpha_{-C}(X)$ by Proposition \ref{PropUCEElementaryProps} (7). The proof for $D_C(z; X)$ runs parallel using Proposition \ref{PropLCEElementaryProps} (6), (7). 
\pend

With the previous result, the first part of the next theorem is immediate. 

\begin{theorem}
\label{ThmExpectileRank}
The expectile ranking functions satisfy:

(1) $D_{-C}$ is lower semicontinuous and quasiconvex, $D_C$ is upper semicontinuous and quasiconcave,

(2) if $A$ is an invertible $d \times d$-matrix and $b \in \R^d$, then $D_{-AC}(Az + b; AX + b) = D_{-C}(z; X)$ and $D_{AC}(Az + b; AX + b) = D_C(z; X)$ for all $X \in  L^1_d$ and all $z \in \R^d$,

(3) $y, z \in \R^d$, $y \leq_C z$ imply $D_{-C}(y; X) \leq D_{-C}(z; X)$ and $D_C(y; X) \leq D_C(z; X)$ for all $X \in L^1_d$,

(4) $X, Y \in L^1_d$, $X \leq_C Y$ a.s. imply $D_{-C}(z; Y) \leq D_{-C}(z; X)$ and $D_C(z; Y) \leq D_C(z; X)$ for all  $z \in \R^d$.
\end{theorem}

{\sc Proof.} (1) Extended real-valued functions are lower (upper) semicontinuous if, and only if, all of their sublevel (superlevel) sets are closed, and they are quasiconvex (quasiconcave) if, and only if, all of their sublevel (superlevel) sets are convex. The result follows from Theorem \ref{ThmExpectileRegionVsRank} and Proposition \ref{PropUCEElementaryProps} (1), Proposition \ref{PropLCEElementaryProps} (1).

(2) One has $D_{-AC}(Az + b; AX + b) = \inf\cb{\alpha \in (0, 1) \mid Az + b \in  E^\alpha_{-AC}(AX+b)}$. The result follows from Proposition \ref{PropUCEElementaryProps} (2) which gives $E^\alpha_{-AC}(AX + b) = AE^\alpha_{-C}(X) + b$.  A parallel argument using Proposition \ref{PropLCEElementaryProps} (2)  works for $D_C$.

(3) The inequality $y \leq_C z$ means $y \in z - C$. If $z \in E^\alpha_{-C}(X)$, then $y \in z - C \subseteq E^\alpha_{-C}(X) - C = E^\alpha_{-C}(X)$ where the last equation is from Proposition \ref{PropUCEElementaryProps} (1). The definition of $D_{-C}(\cdot; X)$ gives the result. The inequality for $D_C(\cdot; X)$ follows similarly with the help of Proposition \ref{PropLCEElementaryProps} (1). 

(4) The claim for $D_{-C}(\cdot; X)$ follows from Proposition \ref{PropUCEElementaryProps} (3), the one for $D_C(\cdot; X)$ from Proposition \ref{PropLCEElementaryProps} (3).\pend

Note that the monotonicity of the expectile ranking functions with respect to the $z$-variable relies on the image space property of $E^\alpha_{-C}$, $E^\beta_C$ whereas the monotonicity with respect to the random variable is a consequence of the corresponding monotonicity of the univariate expectiles.

The expectile depth function $ED$ in \cite[Sec. 4]{CascosOchoa21JMVA} corresponds to the expectile region $ED^\alpha$ and is defined as a supremum since $ED^\alpha$ shrinks with increasing $\alpha$ and thus is a measure of centrality. On the other hand, $E^\alpha_{-C}(X)$ and $E^\beta_C(X)$ "move outward" with decreasing $\alpha$ and increasing $\beta$ in the direction of $-C$ and $C$, respectively. Thus, a low value of $D_{-C}(z; X)$ means that the point $z$ is a "bad" one while a high value of $D_C(z; X)$ means that $z$ is a "good" one if maximality with respect to $\leq_C$ is the goal (vice versa for minimality). Thus, $D_{-C}(z; X)$, $D_C(z; X)$ can be understood as measures for outstandingness in the direction of $-C$ and $C$, respectively. Rather than centrality, the expectile rank functions provide information about the tail behavior of the distribution where "tail" is understood in the direction of $-C$ and $C$, respectively. Therefore, the term "outstandingness" is preferred here over "outlyingness" since the latter does not refer to any order relation for the data points, but is commonly used as the opposite to "central," i.e., being away from the center.

It might be noted that the affine equivariance property \cite[Sec. 4.1 (i)]{CascosOchoa21JMVA} is not convincing since it would imply $ED(0, 0) = ED(z, X)$ for all $z \in \R^d$, $X \in L^1_d$ (choose the zero matrix $A$ and the zero vector $b$) which was certainly not the intention of the authors of \cite{CascosOchoa21JMVA}. Therefore, Theorem \ref{ThmExpectileRank} (2) as well as Proposition \ref{PropUCEElementaryProps} and \ref{PropLCEElementaryProps} assume invertible matrices $A$. In particular, Lemma \ref{LemTransformDualCone} in the Appendix does not work without this assumption.

The following proposition establishes that the downward and the upward expectile rank of a point $z$ generated by $X$ and $-z$ generated by $-X$ sum up to 1 (and hence it is enough to know one of the two as in Proposition \ref{PropLUTransfer} for expectile sets). This shows that the expectile rank functions behave similar to the cumulative distribution function and the survival function of a univariate distribution.

\begin{proposition}
\label{PropSumRanks}
For any $z \in \R^d$ and $X \in L^1_d$, one has
\[ 
D_C(-z; -X) + D_{-C}(z; X) = 1. 
\]
\end{proposition}

{\sc Proof.} By definition and with the help of \eqref{EqLUTransfer} one gets
\[ 
\begin{aligned}
D_{-C}(z; X) & = \inf\cb{\alpha \mid z \in  E^\alpha_{-C}(X) } \\
& = - \sup\cb{-\alpha \mid z \in  E^\alpha_{-C}(X) } = - \sup\cb{-\alpha \mid z \in  - E^{1-\alpha}_C (-X) } \\
& = 1 - \sup\cb{1-\alpha \mid - z \in  E^{1-\alpha}_C (-X) } = 1 - \sup\cb{\beta \mid - z \in  E^{\beta}_C (-X) } \\
& = 1 - D_C (-z; -X)
\end{aligned} 
\]
which is the desired equation. 
\pend

One can define order relations on $\R^d$ using the expectile rank functions as follows:
\begin{align*}
y \leq_{X, -C} z \quad & \Leftrightarrow \quad D_{-C}(y; X) \leq D_{-C}(z; X) \\
y \leq_{X, +C} z \quad & \Leftrightarrow \quad D_C(y; X) \leq D_C(z; X)
\end{align*}
Both $\leq_{X, -C}$ and $\leq_{X, +C}$ are reflexive, transitive and complete since each is defined via one real-valued function. Figure \ref{Cone_Exp_Rank} shows that neither one is antisymmetric in general. Moreover, the monotonicity properties in Theorem \ref{ThmExpectileRank} (3) ensure that $y \leq_C z$ implies $y \leq_{X, -C} z$ as well as $y \leq_{X, +C} z$. However, the example depicted in Figure \ref{Cone_Exp_Rank} shows that both of the above relations do not coincide with $\leq_C$ in general; they both are complete extensions of the vector preoder and thus may serve as (scalarizing) tools in Multi-Criteria Decision Making in the same way as the lower and upper cone distribution function from \cite{HamelKostner18JMVA}: see \cite{HamelKostner23MOD}. Moreover, even the intersection of $\leq_{X, -C}$ and $\leq_{X, +C}$ does not coincide with $\leq_C$ in general: the pair $(y, z)$ in Figure \ref{Cone_Exp_Rank} satisfies $y \leq_{X, -C} z$ as well as $y \leq_{X, +C} z$ (the "orange" $y$ and the "red" $z$ have the same downward and upward expectile rank), but is not comparable with respect to $\leq_{\R^2_+}$.

On the other hand, $(y, z)$  belongs to the symmetric part of $\leq_{X, -C}$ as well as that of $\leq_{X, +C}$ and thus, a decision maker might consider them as equivalent and not prefer one over the other. This motivates the following definition. 

\begin{definition}
\label{DefExpIndifference}
The equivalence relations defined by
\begin{align*}
y \stackrel{X, -C}{\sim} z \quad & \Leftrightarrow \quad D_{-C}(y; X) = D_{-C}(z; X) \\
y \stackrel{X, +C}{\sim} z \quad & \Leftrightarrow \quad D_C(y; X) = D_C(z; X) 
\end{align*}
are called lower and upper expectile rank indifference relation, respectively, with respect to the random variable $X$ and the cone $C$. Their intersection
\[
y \stackrel{X, \pm C}{\sim} z \quad \Leftrightarrow \quad D_{-C}(y; X) = D_{-C}(z; X) \wedge D_C(z; X) = D_C(y; X)
\]
is called expectile rank indifference relation (generated by $X$ and $C$).
\end{definition}

Thus, $y, z \in \R^d$ are lower expectile rank indifferent if $y \leq_{X, -C} z$ and $z \leq_{X, -C} y$ and parallel for $\leq_{X, C}$. Finally, a case is pointed out in which one can draw a conclusion about the comparability of $y, z \in \R^d$ with respect to $\leq_C$ from the knowledge of $\leq_{X, -C}$, $\leq_{X, C}$.

\begin{theorem}
\label{ThmRanksToOrder} 
Let $X \in L^1_d$ and $y, z \in \R^d$ satisfy $D_{-C}(y; X) \leq D_C (z; X)$. Then $y \leq_C z$ if, and only if, $y \leq_{X, -C} z$ and $y \leq_{X, C} z$.
\end{theorem}

{\sc Proof.} "$\Rightarrow$:" This is the monotonicity of the expectile rank functions given in Theorem \ref{ThmExpectileRank} (3) which is true without the additional assumption. "$\Leftarrow$:" Assume
\begin{align*} 
\alpha_1 := D_{-C}(y; X)  & \leq  D_{-C}(z; X) =: \alpha_2, \\
\beta_1: = D_C (y; X) & \leq D_C (z; X) =: \beta_2.
\end{align*}
and hence $\alpha_1 \leq \beta_2$ by assumption. Then $D_{-C}(y; X) \leq \beta_2$ and Theorem \ref{ThmExpectileRegionVsRank} yields $y \in E^{\beta_2}_{-C}(X)$. Thus
\[
\forall w \in C^+ \colon w^\top y \leq e_{\beta_2}(w^\top X).
\]
On the other hand, $D_C (z; X) = \beta_2$ implies $D_C (z; X) \geq \beta_2$, hence $z \in E^{\beta_2}_C(X)$. Thus
\[
\forall w \in C^+ \colon w^\top z \geq e_{\beta_2}(w^\top X).
\]
Together, this yields $w^\top y \leq w^\top z$ for all $w \in C^+$ which means $y \leq_C z$ according to the bipolar theorem.
\pend

Figure \ref{Cone_Exp_Rank} provides a counterexample for the remaining case: one has $D_{-C}(y; X) = D_{-C}(z; X) = 0.35$ and $D_C (y; X) = D_C (z; X) = 0.2$, thus $\alpha_1 > \beta_2$ and the two points $y, z$ are not comparable with respect to $\leq_C$.

This section is concluded with a link between new multivariate stochastic orders and expectile rank functions. 

\begin{definition}
\label{DefExpectileStochOrders}
For two random variables $X, Y \in L^1_d$ two order relations are defined by
\begin{align*}
X \leq_{le} Y \quad \Leftrightarrow \quad & \forall \alpha \in (0, 1) \colon E^\alpha_{-C}(X) \subseteq E^\alpha_{-C}(Y) \\
X \leq_{ue} Y \quad \Leftrightarrow \quad & \forall \alpha \in (0, 1) \colon E^{1-\alpha}_C(X) \supseteq E^{1-\alpha}_C(Y). 
\end{align*}
The relation $\leq_{le}$ is called lower expectile order, the relation $\leq_{ue}$ upper expectile order.
\end{definition}

One may check that in the special case $d = 1$, $C = \R_+$ both $\leq_{le}$ and $\leq_{ue}$ collapse to the expectile order $\leq_e$ defined in \cite[Def. 6]{BelliniKlarMueller18MCAP}. In \cite[Thm. 8]{BelliniKlarMueller18MCAP}, it is shown that $\leq_e$ an be characterized through Omega ratios at least for $L^\infty$ random variables. Instead, the expectile rank functions are used here.

\begin{theorem}
\label{ThmExpectileOrderRank}
One has
\begin{align*}
X \leq_{le} Y \quad \Leftrightarrow \quad & \forall z \in \R^d \colon D_{-C}(z; X) \geq D_{-C}(z; Y) \\
X \leq_{ue} Y \quad \Leftrightarrow \quad & \forall z \in \R^d \colon D_C(z; X) \leq D_C(z; Y). 
\end{align*}
\end{theorem}

{\sc Proof.} First, take $z \in \R^d$ and assume $E^\alpha_{-C}(X) \subseteq E^\alpha_{-C}(Y)$ for all $\alpha \in (0, 1)$. Then, only  two cases are possible:

(i) $E^\alpha_{-C}(Y) = E^\alpha_{-C}(X)$ for all $\alpha \in (0, 1)$ with $z \in E^\alpha_{-C}(Y)$; in this case $D_{-C}(z; X) = D_{-C}(z; Y)$,

(ii) there exist $\alpha \in (0, 1)$ with $z \in E^\alpha_{-C}(Y)\bs E^\alpha_{-C}(X)$; in this case $D_{-C}(z; X) > D_{-C}(z; Y)$.

This proves "$\Rightarrow$" for $\leq_{le}$.

Conversely, fix $\alpha \in (0, 1)$ and take $z \in E^\alpha_{-C}(X)$. Theorem \ref{ThmExpectileRegionVsRank} yields $D_{-C}(z; X) \leq \alpha$. By assumption, $D_{-C}(z; Y) \leq D_{-C}(z; X)$ again by Theorem \ref{ThmExpectileRegionVsRank}, hence $z \in E^\alpha_{-C}(Y)$ as desired.

The proof for $\leq_{ue}$ uses parallel arguments.
\pend

\section{Conclusion and perspective}

In this paper, downward and upward cone expectiles are introduced as set-valued function which capture the tail behavior of a multivariate distribution with respect to an order relation generated by a convex cone instead of its centrality. The cone expectile functions are defined in such a way that they share most of the attractive properties of univariate expectiles such that monotonicity and sublinearity: this is possible since they map into well-defined complete lattices of sets generated by the ordering cone. 

As a result, they also facilitate the contruction of coherent set-valued expectile risk measures which captures the idea of risk compensation much better than the expectile region. The incorporation of makret models as in \cite{HamelHeydeRudloff11MAFE} would be the natural next step.

 The dual information shared by both cone expectiles in the present paper and the expectile region from \cite{DiGiorgiMcNeil16CSDA, CascosOchoa21JMVA} opens dual options for computing these sets. The design of efficient algorithms along with convergence properties as well as questions of the asymptotic behaviour of cone expectiles remains a task.

The new downward and upward expectile rank functions correspond to the lower and upper expectile in the same way as the univariate cumulative distribution function corresponds to the quantile function. They admit to quantify the "outstandingness" of a (data) point with respect to a multivariate distribution.

The new upper and lower cone expectiles as well as the corresponding rank functions can be applied to the analysis of ordered data and especially to the ranking problem in Multi-Criteria Decision Making.

\bibliographystyle{plain}

\renewcommand{\thelemma}{\Alph{subsection}.\arabic{lemma}}
\renewcommand{\thesubsection}{\Alph{subsection}.}

\section*{Appendix}

\subsection{Lemmas}

\begin{lemma}
\label{LemTransformDualCone}
Let $C \subseteq \R^d$ be a nontrivial closed convex cone and $A$ an invertible $d \times d$-matrix. Then $AC$ is a closed convex cone with $(AC)^+ = (A^\top )^{-1}C^+$. 
\end{lemma}

{\sc Proof.} It is straighforward to check that $AC$ is a closed convex cone. One has
\begin{align*}
v \in (AC)^+ & \Leftrightarrow \forall z \in C \colon v^\top  (Az) = (A^\top  v)^\top  z \geq 0 \\
	& \Leftrightarrow  A^\top  v \in C^+ \Leftrightarrow v \in (A^\top )^{-1}C^+
\end{align*}
which completes the proof. \pend

\begin{lemma}
\label{LemHalfSpaceSum}
If $w \in \R^d\bs\{0\}$ and $r, s \in \R$, then
\[
\cb{z \in \R^d \mid w^\top  z \leq r + s} = \cb{z \in \R^d \mid w^\top  z \leq r} + \cb{z \in \R^d \mid w^\top  z \leq s}
\]
where the addition is understood as element-wise (Minkowski) addition of subsets of $\R^d$.
\end{lemma}

{\sc Proof.} Since $w \neq 0$, there are $x(r), x(s) \in \R^d$ with $w^\top  x(r) = r$ and $w^\top  x(s) = s$. Then
\[
\cb{z \in \R^d \mid w^\top  z \leq r} = H^+(-w) + \{x(r)\}
\]
and the parallel equation is true for $s$. Indeed, if $z = x + x(r)$ with $-w^\top  x \geq 0$, then $w^\top  z = w^\top  (x + x(r)) \leq r$, i.e., $z$ belongs to the left hand side. If $w^\top  z \leq r$, then $z = z - x(r) + x(r)$ with $w^\top (z - x(r)) \leq r - r = 0$, i.e., $z - x(r) \in H^+(-w)$, and therefore, $z$ belongs to the right hand side. \pend

\subsection{Proof of Theorem \ref{PolyhedralConeExp}} 
\label{Prof_PolConeExp}

The proof is given for the downward cone expectile $E^{\alpha}_{-C}(X)$. The proof for $E^{1-\alpha}_C(X)$ follows similar lines.
Define the set
\[
\widehat E^{\alpha}_{-C}(X) = 
\bigcap_{m \in \cb{1, \ldots, M}} \cb{z \in \R^{d} \mid (w^m)^\top z \leq e_{\alpha}((w^m)^\top X)}.
\]
Clearly, $\widehat E^{\alpha}_{-C}(X) \supseteq E^{\alpha}_{- C}(X)$, so the opposite inclusion needs to be verified. Since the set $\{w^1, \ldots, w^M\}$ generates $C^+$, every $w \in C^+$ can be represented as $w = \sum_{m = 1}^M s_m w^m$ with numbers $s_1, \ldots, s_M \geq 0$. Since $0< \alpha \leq \frac{1}{2}$, the univariate expectiles $e_{\alpha}$ are super-additive and positively homogeneous which yields
\[
e_{\alpha}(w^\top X) \geq \sum_{i=1}^{M} s_m e_{\alpha}((w^m)^\top X).
\]
In turn, this leads to 
\begin{equation}
\label{EqInclusion1}
\cb{z \in \R^d \mid w^\top z \leq e_{\alpha}(w^\top X)} \supseteq
	\cb{z \in \R^d \mid w^\top z \leq \sum_{m=1}^{M} s_m e_{\alpha}((w^m)^\top X)}
\end{equation}
If $y \in \widehat E^{\alpha}_{-C}(X)$ one has $s_m (w^m)^\top y \leq s_m e_{\alpha}((w^m)^\top X)$ for all $m \in \{1, \ldots, M\}$.
Therefore,
\[
w^\top y = \left[\sum_{m = 1}^{M} s_m w^m \right]^\top y \leq \sum_{m=1}^{M} s_m e_{\alpha}((w^m)^\top X),
\]
which gives $y \in \cb{z \in \R^d \mid w^\top z \leq \sum_{m=1}^{M} s_me_{\alpha}((w^m)^\top X)}$ and consequently
\begin{equation}
\label{EqInclusion2}
\widehat E^{\alpha}_{-C}(X) \subseteq \cb{z \in \R^d \mid w^\top z \leq \sum_{m=1}^{M} s_me_{\alpha}((w^m)^\top X)} 
\end{equation}
The two inclusions \eqref{EqInclusion1} and \eqref{EqInclusion2} give
\begin{equation}
\label{EqInclusion3}
\widehat E^{\alpha}_{-C}(X) \subseteq \cb{z \in \R^d \mid w^\top z \leq e_{\alpha}(w^\top X)}. 
\end{equation}
Since \eqref{EqInclusion3} holds for all $w \in C^+$, one can take the intersection over $w \in C^+$ and obtains $\widehat E^{\alpha}_{-C}(X) \subseteq E^{\alpha}_{- C}(X)$ which completes the proof. \pend

\end{document}